\documentclass[12pt,draftcls,onecolumn]{IEEEtran}
\usepackage{amssymb}
\usepackage{graphicx,subfigure}
\usepackage{latexsym,amsmath,amsthm,amssymb,amsfonts}
\usepackage{mathrsfs}
\usepackage{amscd}
\usepackage[ansinew]{inputenc}
\usepackage{color}
\usepackage{cite}
\usepackage[all]{xy}

\def\QED{\mbox{\rule[0pt]{1.5ex}{1.5ex}}}
\def\endproof{\hspace*{\fill}~\QED\par\endtrivlist\unskip}

\newtheorem{theorem}{Theorem}

\newtheorem{remark}{Remark}
\newtheorem{lemma}{Lemma}
\newtheorem{assumption}{Assumption}

\begin{document}

\title {Distributed Adaptive Gradient Optimization Algorithm
\thanks{This work was supported by the National Science Foundation under Grant ECCS-1307678 and ECCS-1611423, and the National Natural Science Foundation of China (61203080,61573082,61528301).}
}
\author{\mbox{Peng Lin}\thanks{Peng Lin is with the School of Information Science and Engineering, Central South University, Changsha 410083, China. Wei Ren is with the Department
of Electrical and Computer Engineering, University of California, Riverside, CA92521, USA.
E-mail: lin$\_$peng0103@sohu.com,
ren@ee.ucr.edu.} and Wei Ren}

\maketitle

\begin{abstract}
In this paper, a distributed optimization problem with general differentiable convex objective functions is
studied for single-integrator and double-integrator multi-agent systems. Two distributed adaptive optimization algorithm is introduced which uses the
relative information to construct the gain of the interaction term. 
The analysis is performed based on the Lyapunov functions, the analysis of the system solution and the convexity of the local objective functions. It is shown that if the gradients of the convex objective functions are continuous, the team convex objective function can be minimized as time evolves for both single-integrator and double-integrator multi-agent systems.
Numerical examples are
included to show the obtained theoretical results.
\vspace{0.1cm}

\noindent{\bf Keywords}: Optimization, Consensus, Distributed Adaptive algorithm
\end{abstract}

\section{Introduction}
As an important branch of distributed control, distributed optimization has
attracted
more and more attention from the control community
\cite{angelia,angelia1,Zhu,Johansson,linren3,lu,lup,shi,liu,Kvaternik,Elia,cotes,Cortes3,linren1,linren2,zhao}. The aim is to use a distributed approach to minimize a team optimization function composed of a sum of local objective functions where each local objective function is known to only one agent. In the past few years, researchers have obtained many results about distributed optimization problems from different perspectives. For example, based on gradient descent method, articles  \cite{angelia,angelia1,shi,liu,linren3} studied distributed optimization problems with and without state constraints, while  by introducing a dynamic integrator, articles \cite{Elia,cotes,Cortes3} investigated distributed optimization problems for general strongly connected balanced directed graphs. Recently, some researchers turned their attention to try to solve the distributed optimization problem from a view point of nonsmooth approaches. For example, article \cite{linren1} proposed several algorithms using nonsmooth functions to solve the distributed optimization problem with the consideration of finite-time consensus optimization convergence. Also, articles \cite{linren3,zhao} introduced  adaptive algorithms using nonsmooth functions to solve the distributed optimization problem for general differentiable convex functions or general linear multi-agent systems. However, in \cite{linren1,linren3,zhao}, it is required that the gradients or subgradients of the local objective functions be bounded or a period of the previous information should be used for each agent.

To this end, we will continue the work of \cite{linren3} to study the distributed optimization problem for general differentiable objective function using nonsmooth functions.
Two distributed adaptive optimization algorithm is introduced which uses the relative information to construct the gain of the interaction term. The analysis is performed based on the Lyapunov functions, the analysis of the system solution and the convexity of the local objective functions. It is shown that if the gradients of the convex objective functions are continuous, the team convex objective function can be minimized as time evolves for both single-integrator and double-integrator multi-agent systems.

 Notations. $\mathbb{R}^m$  denotes the set of all $m$ dimensional real column
vectors; $\mathcal{I}$ denotes the index set
$\{1,\ldots,n\}$; $s_i$ denotes the $i$th component of the vector $s$;
$s^T$ denotes the transpose of the vector $s$; $||s||$ denotes the Euclidean norm of the vector $s$; {{$\frac{\mathrm{d}}{\mathrm{d} s}$}}
denotes the differential operator with respect to $s$; {$\nabla
f(s)$ denotes the gradient of the function $f(s)$ at $s$;}
  $\mathrm{sgn}(s)$ denotes a component-wise sign function of $s$; and {$P_{X}(s)$} denotes the projection of the vector $s$ onto the closed convex set $X$, i.e., {$P_{X}(s)=\mathrm{arg} \min\limits_{\bar{s}\in
X}\|s-\bar{s}\|$}.

\section{Preliminaries}
In this section, we introduce preliminary results
about graph theory and convex functions (see \cite{s10, boyd, linren3}).

Consider a multi-agent system consisting of $n$ agents. Each agent is regarded as a node in an undirected graph $\mathcal{G}(\mathcal{V},\mathcal{E},\mathcal{A})$ of order $n$ where $\mathcal{V}=\{1,\cdots,n\}$ is
the set of nodes, $\mathcal{E}\subseteq\mathcal{V}\times
\mathcal{V}$ is the set of edges, and $\mathcal{A}=[a_{ij}]\in
\mathbb{R}^{n\times n}$ is the weighted adjacency matrix. An edge of
$(i,j) \in \mathcal{E}$ denotes that agents  $i$ and $j$  can obtain
information from each other. The weighted adjacency matrix $\mathcal{A}$ is
defined as $a_{ii}=0$ and $a_{ij}=a_{ji}\neq 0$ if $(i,j) \in
\mathcal{E}$ and $a_{ij}=0$ otherwise. The set of neighbors of node
$i$ is denoted by $N_i=\{j\in \mathcal{V}:(i,j)\in \mathcal{E}\}.$ The Laplacian of the graph $\mathcal{G}$, denoted by $L$,
is defined as $\lfloor {L}\rfloor_{ii} =\sum_{j=1}^na_{ij}$ and $\lfloor {L}\rfloor_{ij}=-a_{ij}$ for all $i\neq j$. A
path is a sequence of edges of the form
$(i_1,i_2),(i_2,i_3),\cdots$, where $i_j \in \mathcal{V}$. The graph $\mathcal{G}$ is  connected, if there
is a path from every node to every other node.
\begin{lemma}\label{lemma1}
    {\rm\cite{s10}
    $~$If the graph $\mathcal{G}$ is connected, then its Laplacian $L$ has a simple eigenvalue at $0$ with associated eigenvector
     $\textbf{1}$ and all its other $n-1$ eigenvalues are positive and real.
     }
    \end{lemma}

\begin{lemma}\label{lemma35}{\rm\cite{boyd} Let $f_0(s): \mathbb{R}^m\rightarrow\mathbb{R}$ be a differentiable convex function. $f_0(s)$ is minimized if and only if $\nabla f_0(s)=0$.}\end{lemma}

\begin{lemma}\label{lemmallin}{\rm \cite{linren3} Under Assumption \ref{ass2},  all $X_i$ and $X$ are nonempty closed bounded convex sets for all $i$.}\end{lemma}

\section{Distributed Optimization Problem}
Suppose
that each agent has the following dynamics
\begin{equation} \label{eq:single-int}
\dot{x}_i(t)=u_i(t), \quad i\in \mathcal{I},
\end{equation}
where $x_i(t)\in \mathbb{R}^m$ is the state of agent $i$, and $u_i(t) \in
\mathbb{R}^m$ is the control input of agent $i$.
Our objective is to use only local information to
design $u_i(t)$ for all agents to cooperatively solve
the following optimization problem
\begin{eqnarray}\label{gel16}\begin{array}{lll}\mathrm{minimize}~~\sum_{i=1}^nf_i(x_i)\\
\mathrm{subject~to}~~x_i=x_j\in
\mathbb{R}^m.\end{array}\end{eqnarray}

\begin{assumption}\label{ass2}{\rm \cite{linren3} Each set $X_i\triangleq{\big\{}s{\big |}\nabla f_i(s)=0{\big\}}$ is nonempty and bounded.
%
}\end{assumption}

\begin{assumption}\label{ass24}{\rm \cite{linren3} The length of the time interval between any two contiguous switching times is no smaller than a given constant, denoted by $d_w$. }\end{assumption}

\section{Main Results}
\subsection{Single-Integrator Multi-Agent Systems}
In this section, we  design a distributed adaptive algorithm for (\ref{eq:single-int}) to solve the optimization problem (\ref{gel16}) for general convex local objective functions.
The algorithm is given by
\begin{eqnarray}\label{gel2}\begin{array}{lll}
u_i(t)=\sum_{j\in
N_i(t)}\frac{q_{ij}(t)[x_j(t)-x_i(t)]}{\|x_j(t)-x_i(t)\|}-\nabla
f_i(x_i(t)),\\
\dot{q}_{ij}(t)=\left\{\begin{array}{lll}\mathrm{sgn}(\|x_j(t)-x_i(t)\|),\mbox{~if~}(i,j)\in \mathcal{G}(t),\\0,\mbox{~otherwise,~}\end{array}\right.\\q_{ij}(0)=q_{ji}(0)=0,\end{array}
\end{eqnarray}
for all $i$. In (\ref{gel2}), the role of the term, $\sum_{j\in
N_i(t)}\frac{q_{ij}(t)[x_j(t)-x_i(t)]}{\|x_j(t)-x_i(t)\|}$, is to make all agents converge to a consensus point, while the second term,
$-\nabla
f_i(x_i(t))$, is the negative gradient of $f_i(x_i(t))$ which is used to minimize $f_i(x_i(t))$.

\begin{remark}{\rm As algorithm (\ref{gel2}) uses
 the sign functions that is nonsmooth, the system
(\ref{eq:single-int}) using (\ref{gel2}) would be discussed in the Filippov sense \cite{Fili}.}\end{remark}



\begin{theorem}\label{theorem1}{\rm Suppose
that {{the graph $\mathcal{G}(t)$ is undirected and connected for all $t$, $\nabla f_i(s)$ is continuous with respect to $s$ for all $i$ and Assumptions \ref{ass2} and \ref{ass24} hold}}.
 {{For system
(\ref{eq:single-int}) with algorithm (\ref{gel2}),}} all agents reach a
consensus in finite time and minimize the
team objective function (\ref{gel16}) as $t\rightarrow+\infty$.
}
\end{theorem}
\noindent{\textbf{Proof:}} First, we prove that all $x_i(t)$ are bounded for all $t$. Under Assumption \ref{ass2}, from Lemma \ref{lemmallin}, we have that all $X_i$ and $X$ are
 nonempty closed bounded convex sets for all $i$. It is clear that $x_i(0)\in Y$,
$X\subset Y$ and $X_i\subset Y$ for all $i$ and some closed bounded set $Y$. Let $Y$ be sufficiently large for any $z\in X$and all $z_j\in X_j$ such that $f_i(x_i(t))-f_{i}(z)\geq \sum_{j=1,j\neq i}^n[f_j(z)-f_j(z_j)]$ for all $i$.
Since $z\in X$, from the convexity of the function $f_i(x_i(t))$, we have $\nabla f_i(x_i(t))^T(z-x_i(t))\leq f_i(z)-f_i(x_i(t))$.

Construct a Lyapunov function candidate as
$V(t)=\frac{1}{2}\sum_{i=1}^n\|x_i(t)-z\|^2$ for some $z\in X$.
Calculating $\dot{V}(t)$ along the solutions of system (\ref{eq:single-int}) with (\ref{gel2}), we have
{\begin{eqnarray}\label{eq551}\begin{array}{lll}&&\dot{V}(t)\\&=&\sum_{i=1}^n[x_i(t)-z]^T\\&\times&
{\big[}\sum_{j\in N_i(t)}\frac{q_{ij}(t)[x_j(t)-x_i(t)]}{\|x_j(t)-x_i(t)\|}-\nabla
f_i(x_i(t)){\big]}
\end{array}\end{eqnarray}}
Since the graph $\mathcal{G}(t)$ is undirected, it follows that
{\begin{eqnarray}\label{eq88pp}\begin{array}{lll}&&\sum_{i=1}^n[x_i(t)-z]^T\\&\times&\sum_{j\in
N_i(t)}\frac{q_{ij}(t)[x_j(t)-x_i(t)]}{\|x_j(t)-x_i(t)\|}\\&=&\sum_{i=1}^n\sum_{j\in
N_i(t)}q_{ij}(t)[x_i(t)-z]^T\\&\times&\frac{x_j(t)-x_i(t)}{\|x_j(t)-x_i(t)\|}\\&=&\sum_{i=1}^n\sum_{j\in
N_i(t)}{\big{\{}}\frac{q_{ij}(t)}{2}\\&\times&[x_i(t)-z]^T\frac{x_j(t)-x_i(t)}{\|x_j(t)-x_i(t)\|}\\
&+&\frac{q_{ij}(t)}{2}[x_j(t)-z]^T\frac{x_i(t)-x_j(t)}{\|x_j(t)-x_i(t)\|}{\big{\}}}\\
&=&\sum_{i=1}^n\sum_{j\in
N_i(t)}\frac{q_{ij}(t)}{2}[x_i(t)-z\\&-&x_j(t)+z]^T\frac{x_j(t)-x_i(t)}{\|x_j(t)-x_i(t)\|}\\
&=&\sum_{i=1}^n\sum_{j\in
N_i(t)}\frac{q_{ij}(t)}{2}[x_i(t)-x_j(t)]^T\\&\times&\frac{x_j(t)-x_i(t)}{\|x_j(t)-x_i(t)\|}\\
&\leq&0.\end{array}\end{eqnarray}}
From the convexity of the function $f_i(x_i(t))$, we have $\nabla f_i(x_i(t))^T(z-x_i(t))\leq f_i(z)-f_i(x_i(t))$.
 It follows that $\dot{V}(t)\leq -\sum_{i=1}^n[f_i(x_i(t))-f_i(z)]$. If $x_{i_0}(t)\notin Y$ for some $i_0$, we have $f_{i_0}(x_{i_0}(t))-f_{i_0}(z)\geq \sum_{j=1,j\neq {i_0}}^n[f_j(z)-f_j(z_j)]$ for all $z_j\in X_j$ and hence $\dot{V}(t)\leq -[f_{i_0}(x_i(t))-f_{i_0}(z)]+\sum_{j=1,j\neq i_0}^n[f_i(z)-f_i(z)]\leq 0$. This
implies that all $x_i(t)$ remain in $Y$. Note that each $\nabla f_i(x_i(t))$ is continuous with respect to $x_i(t)$ for all $i$, $X\subset Y$ and $Y$ is bounded. Thus, $\max\{\|x_i(t)\|,\|\nabla f_i(t)\|\}<\rho$ for all $i$ and some constant $\rho>0$.

Next, we prove that all agents reach a consensus as $t\rightarrow+\infty$. Let $0<t_{k1}<t_{k2}<t_{k+1,1}<t_{k+1,2}$ denote the contiguous switching times for all $k\in \{1,2,\cdots\}$ such that $x_i(t)\neq x_j(t)$ for some two integers $i,j\in \mathcal{I}$ and all $t\in[t_{k1},t_{k2})$ and $x_i(t)=x_j(t)$ for all $i,j\in \mathcal{I}$ and all $t\in[t_{k2},t_{k+1,1})$.
Suppose that consensus is not reached as $t\rightarrow+\infty$ and $\sum_{k=1}^{+\infty}(t_{k2}-t_{k1})<+\infty$. It is clear that
$\lim_{k\rightarrow+\infty}(t_{k2}-t_{k1})=0$. Moreover, from the dynamics of $q_{ij}(t)$, we have that $q_{ij}(t)<\rho_q$ for some constant $\rho_q>0$. Since each $\nabla f_i(s)$ is bounded and $x_i(t_{k1}^-)=x_j(t_{k1}^-)$ for all $i,j$ and all $k>1$, where $t_{k1}^-$ denotes the time just before $t_{k1}$, $\|u_i(t)\|$ is bounded for all $t\in[t_{k1},t_{k2})$ and hence $0\leq\lim_{k\rightarrow+\infty}\max_{t\in[t_{k1},t_{k2})}\|x_i(t)-x_j(t)\|\leq \lim_{k\rightarrow+\infty}\int_{t_{k1}}^{t_{k2}}(\|u_i(s)\|+\|u_j(s)\|)\mathrm{d}s=0$ for all $i,j$. That is, consensus is reached as $t\rightarrow+\infty$, which yields a contradiction. Suppose that $\sum_{k=1}^{+\infty}(t_{k2}-t_{k1})=+\infty$.
Similar to the proof of Theorem 2 in \cite{linren3}, it can be proved that all agents reach a consensus in finite time.

Summarizing the above analysis, consensus can be reached as $t\rightarrow+\infty$. Let $x^*(t)=\frac{1}{n}\sum_{i=1}^nx_i(t)$. Note that each $\nabla f_i(x_i(t))$ is continuous with respect to $x_i(t)$. There is a constant $T>0$ for any $\epsilon>0$ such that $\|x_i(t)-x^*(t)\|<\epsilon$ and $\|\nabla f_i(x^*(t))-\nabla f_i(x_i(t))\|<\epsilon$ for all $t>T$.
Recall that $\|x_i(t)\|<\rho$.
Consider the Lyapunov function candidate $V_1(t)=\frac{1}{2}\|x^*(t)-P_X(x^*(t))\|^2$ for $t>T$. Calculating $\dot{V}_1(t)$, we have
\begin{eqnarray*}\dot{V}_1(t)=-[x^*(t)-P_X(x^*(t))]^T\frac{1}{n}\sum_{i=1}^n\nabla f_i(x_i(t))\\
=-[x^*(t)-P_X(x^*(t))]^T[\frac{1}{n}\sum_{i=1}^n\nabla f_i(x^*(t))+\epsilon]\\
\leq -[\frac{1}{n}\sum_{i=1}^n f_i(x^*(t))-\frac{1}{n}\sum_{i=1}^n f_i(P_X(x^*(t)))]+2\rho\epsilon\end{eqnarray*}
Note that when $\frac{1}{n}\sum_{i=1}^n f_i(x^*(t))-\frac{1}{n}\sum_{i=1}^n f_i(P_X(x^*(t)))\geq 4\rho\epsilon$, $\dot{V}_1(t)\leq -2\rho\epsilon$. It follows that there exists a constant $T_1>T$ such that $\frac{1}{n}\sum_{i=1}^n f_i(x^*(t))-\frac{1}{n}\sum_{i=1}^n f_i(P_X(x^*(t)))< 4\rho\epsilon$ for $t>T_1$. Since $\epsilon$ can be arbitrarily small, it follows that $\lim_{t\rightarrow+\infty}[\frac{1}{n}\sum_{i=1}^n f_i(x^*(t))-\frac{1}{n}\sum_{i=1}^n f_i(P_X(x^*(t)))]=0$. It follows from Lemma \ref{lemma35} that the
team objective function (\ref{gel16}) is minimized as $t\rightarrow+\infty$.\endproof

\begin{remark}{\rm In \cite{linren3}, a distributed algorithm was proposed to solve the optimization problem. However, it is required that a period of the previous information should be used for each agent. In contrast to \cite{linren3}, in this paper, the previous information is not used and the current information is sufficient for the proposed algorithm to make all agents minimize the team objective function as time evolves. }\end{remark}

\subsection{Double-Integrator Multi-agent Systems}
In this part, our goal is to extend the results in Subsection A to second-order multi-agent systems with the following dynamics
\begin{eqnarray}\label{equationsec1}\begin{array}{lll}
\dot{x}_i(t)&=&v_i(t)\\
\dot{v}_i(t)&=&u_i(t),
\end{array}\end{eqnarray}
where $x_i(t)\in \mathbb{R}^m$ and $v_i(t)\in \mathbb{R}^m$ are
the position and velocity states of agent $i$ and $u_i(t)\in \mathbb{R}^m$ is
the control input. To solve the distributed optimization problem, we use the following algorithm
\begin{eqnarray}\label{equationsec2}\begin{array}{lll}
u_i(t)=-pv_i(t)+\sum_{j\in N_i(t)}\frac{q_{ij}(t)[x_j(t)+\frac{2}{p}v_j(t)-x_i(t)-\frac{2}{p}v_i(t)]}{\|x_j(t)+\frac{2}{p}v_j(t)-x_i(t)-\frac{2}{p}v_i(t)\|}-\nabla f_i(x_i(t)+\frac{2}{p}v_i(t)),\\
\dot{q}_{ij}(t)=\left\{\begin{array}{lll}\mathrm{sgn}(\|x_j(t)-x_i(t)\|),\mbox{~if~}(i,j)\in \mathcal{G}(t),\\0,\mbox{~otherwise,~}\end{array}\right.\\q_{ij}(0)=q_{ji}(0)=0,
\end{array}\end{eqnarray}
where $p>0$ is the feedback damping gain of the agents.

Let $\bar{v}_i(t)=x_i(t)+\frac{2v_i(t)}{p}$. The system (\ref{equationsec1}) with (\ref{equationsec2}) can be written as
\begin{eqnarray}\label{equationsec3}\begin{array}{lll}
\dot{x}_i(t)&=&\frac{p}{2}\bar{v}_i(t)-\frac{p}{2}x_i(t)\\
\dot{\bar{v}}_i(t)&=&-\frac{p}{2}\bar{v}_i(t)+\frac{p}{2}x_i(t)+\frac{2}{p}\sum_{j\in N_i(t)}\frac{q_{ij}(t)[\bar{v}_j(t)-\bar{v}_i(t)]}{\|\bar{v}_j(t)-\bar{v}_i(t)\|}-\frac{2}{p}\nabla f_i(\bar{v}_i(t)).
\end{array}\end{eqnarray}
For convenience of expression, we assume $m=1$ in the proof of the following theorem.
\begin{theorem}\label{theorem2}{\rm Suppose
that {{the graph $\mathcal{G}(t)$ is undirected and connected for all $t$, $\nabla f_i(s)$ is continuous with respect to $s$ for all $i$ and Assumptions \ref{ass2} and \ref{ass24} hold}}.
 {{For system
(\ref{eq:single-int}) with algorithm (\ref{gel2}),}} all agents reach a
consensus in finite time and minimize the
team objective function (\ref{gel16}) as $t\rightarrow+\infty$.
}
\end{theorem}
\noindent{\textbf{Proof:}} 
Construct a Lyapunov function candidate as
$V(t)=\frac{1}{2}\sum_{i=1}^n\|x_i(t)-s\|^2+\frac{1}{2}\sum_{i=1}^n\|\bar{v}_i(t)-s\|^2$ for some $s\in X$. 
Let $z(t)=[x_1(t)^T,\bar{v}_1(t)^T,\cdots,x_n(t)^T,\bar{v}_n(t)^T]^T$, $A=\begin{bmatrix}\frac{p}{2}&-\frac{p}{2}\\-\frac{p}{2}&\frac{p}{2}\end{bmatrix}$, 
$B=\begin{bmatrix}0&0\\0&\frac{2}{p}\end{bmatrix}$ and $\Phi(t)$ be a matrix with each entry $[\Phi(t)]_{ij}=\left\{\begin{array}{lll}-\sum_{k=1,k\neq i}^n[\Phi(t)]_{ik}, \mbox{if}~i=j,\\-\frac{q_{ij}(t)}{2\|x_j(t)-x_i(t)\|}, \mbox{if}~i\neq j~\mbox{and}~(i,j)\in \mathcal{E}(\mathcal{G}(t))\\0,\mbox{otherwise}.\end{array}\right.$
Regarding $A$ and $\Phi(t)$ as the Laplacians of some certain undirected graphs, it follows from Lemma \ref{lemma1} that $-z(t)^T(I_n\otimes A)z(t)\leq 0$ and $-z(t)^T[\Phi(t)\otimes B]z(t)\leq 0$.

Calculating $\dot{V}(t)$, we have
\begin{eqnarray*}\begin{array}{lll}\dot{V}(t)&=-z(t)^T(I_n\otimes A)z(t)-z(t)^T[\Phi(t)\otimes B]z(t)\\&-\sum_{i=1}^n\frac{2}{p}(\bar{v}_i(t)-s)^T \nabla
f_i(\bar{v}_i(t))\\&\leq -\sum_{i=1}^n\frac{2}{p}\|v_i(t)\|^2-\frac{2}{p}\sum_{i=1}^n[f_i(\bar{v}_i(t))-f_i(s)]\\&-z(t)^T[\Phi(t)\otimes B]z(t)\\&\leq-\sum_{i=1}^n\frac{2}{p}\|v_i(t)\|^2-\frac{2}{p}\sum_{i=1}^n[f_i(\bar{v}_i(t))-f_i(s)]\\&-z(t)^T[\Phi(t)\otimes B]z(t),\end{array}\end{eqnarray*}
where the first inequality uses the convexity of $f_i(\cdot)$.
Then by a similar approach to the proof of Theorem 1, it can be proved that all $x_i(t)$ and $\bar{v}_i(t)$ remain in a bounded closed convex set, denoted by $Y$, for all $t$ such that $X\in Y$, $X_i\in Y$ and $x_i(0)\in Y$ for all $i$.
Note that each $\nabla f_i(\bar{v}_i(t))$ is continuous with respect to $\bar{v}_i(t)$. Thus, $\max\{\|x_i(t)\|,\|\bar{v}_i(t)\|,\|\nabla f_i(\bar{v}_i(t))\|\}<\rho$ for all $i$ and some constant $\rho>0$.

Next, we prove that all agents reach a consensus as $t\rightarrow+\infty$. Let $0<t_{k1}<t_{k2}<t_{k+1,1}<t_{k+1,2}$ denote the contiguous switching times for all $k\in \{1,2,\cdots\}$ such that $x_i(t)\neq x_j(t)$ for some two integers $i,j\in \mathcal{I}$ and all $t\in[t_{k1},t_{k2})$ and $x_i(t)=x_j(t)$ for all $i,j\in \mathcal{I}$ and all $t\in[t_{k2},t_{k+1,1})$.
Suppose that consensus is not reached as $t\rightarrow+\infty$ and $\sum_{k=1}^{+\infty}(t_{k2}-t_{k1})<+\infty$. It is clear that
$\lim_{k\rightarrow+\infty}(t_{k2}-t_{k1})=0$. Moreover, from the dynamics of $q_{ij}(t)$, we have that $q_{ij}(t)<\rho_q$ for some constant $\rho_q>0$. Note that $\max\{\|x_i(t)\|,\|\bar{v}_i(t)\|\}<\rho$ for all $i$ and $x_i(t_{k1}^-)=x_j(t_{k1}^-)$ for all $i,j$ and all $k>1$, where $t_{k1}^-$ denotes the time just before $t_{k1}$. Hence $0\leq\lim_{k\rightarrow+\infty}\max_{t\in[t_{k1},t_{k2})}\|x_i(t)-x_j(t)\|=0$ for all $i,j$. That is, $\lim_{t\rightarrow+\infty}[x_i(t)-x_j(t)]=0$ for all $i,j$.
 Since
$x_i(t)=x_j(t)$ for all $i,j\in \mathcal{I}$ and all $t\in[t_{k2},t_{k+1,1})$, it follows from the dynamics of each agent that $v_i(t)=v_j(t)$ for all $i,j\in \mathcal{I}$ and all $t\in(t_{k2},t_{k+1,1})$. Since $q_{ij}(t)<\rho_q$ and $\max\{\|x_i(t)\|,\|\bar{v}_i(t)\|,\|\nabla f_i(\bar{v}_i(t))\|\}<\rho$ for all $i$, it follows that each $u_i(t)$ is bounded for all $i$. Hence
$0\leq\lim_{k\rightarrow+\infty}\max_{t\in[t_{k1},t_{k2})}\|v_i(t)-v_j(t)\|\leq \lim_{k\rightarrow+\infty}\int_{t_{k1}}^{t_{k2}}(\|u_i(s)\|+\|u_j(s)\|)\mathrm{d}s=0$ for all $i,j$. 

Recall that $x_i(t)=x_j(t)$  for all $i,j$ and all $t\in (t_{k2},t_{k+1,1})$. Clearly, $\dot{V}(t)\leq-\sum_{i=1}^n\frac{2}{p}\|v_i(t)\|^2$ for all $t\in (t_{k2},t_{k+1,1})$.
Since $q_{ij}(t)<\rho_q$ and $\max\{\|x_i(t)\|,\|\bar{v}_i(t)\|,\|\nabla f_i(\bar{v}_i(t))\|\}<\rho$ for all $i$, $\dot{V}(t)$ is bounded for all $t$. Since $\sum_{k=1}^{+\infty}(t_{k2}-t_{k1})<+\infty$, $\sum_{k=1}^{+\infty}\int_{t_{k1}}^{t_{k2}}\dot{V}(s)\mathrm{d}s$ is bounded. Thus, $V(t)$ is bounded for all $t$. Thus, $\sum_{k=1}^{+\infty}\int_{t_{k2}}^{t_{k+1,1}}\dot{V}(s)\mathrm{d}s\leq -\sum_{k=1}^{+\infty}\sum_{i=1}^n\int_{t_{k2}}^{t_{k+1,1}}\frac{2}{p}\|v_i(s)\|^2\mathrm{d}s$ is also bounded. This means that $\lim_{k\rightarrow+\infty}\max_{t\in[t_{k2},t_{k+1,1}]}\|v_i(t)\|=0$. By a similar approach to prove that $\lim_{t\rightarrow+\infty}[x_i(t)-x_j(t)]=0$ for all $i,j$, using the continuity of $v_i(t)$, it can be proved that
$\lim_{t\rightarrow+\infty}v_i(t)=0$ for all $i$. It follows from the definition of $\bar{v}_i(t)$ that
$\lim_{t\rightarrow+\infty}[\bar{v}_i(t)-{x}_i(t)]=\lim_{t\rightarrow+\infty}[\bar{v}_i(t)-\bar{v}_j(t)]=0$ for all $i,j$.


Suppose that $\sum_{k=1}^{+\infty}(t_{k2}-t_{k1})=+\infty$. Then from the dynamics of $q_{ij}(t)$, there must exist a pair of agents, denoted by $i_0\neq j_0$, such that $\lim_{t\rightarrow+\infty}q_{i_0j_0}(t)=+\infty$. In the following, we prove that there exist a pair of agents, denoted by $i_1\neq j_1$, such that $(i_1,j_1)\notin\{(i_0,j_0),(j_0,i_0)\}$, $i_1\in \{i_0,j_0\}$ and $\lim_{t\rightarrow+\infty}q_{i_1j_1}=+\infty$. If this is not true, we have $q_{ii_0}(t)<\gamma_q$ and $q_{ij_0}(t)<\gamma_q$ for some constant $\gamma_q>\max_i\{\rho, p\rho\}$, all $t$ and all $i\in \cup_{s\in [0,+\infty)} [N_{i_0}(s)\cup N_{j_0}(s)]$ with $i\neq i_0$ and $i\neq j_0$. Since $\lim_{t\rightarrow+\infty}q_{i_0j_0}(t)=+\infty$, there exists a sufficiently large constant $T_0>0$ for any $\gamma_0>16nm\gamma_q$ such that $q_{i_0j_0}(t)>\gamma_0$ for all $t>T_0$. 
By simple calculations based on (\ref{equationsec3}), when $(i_0,j_0)\in \mathcal{G}(t)$ and $\|x_{i_0}(t)-x_{j_0}(t)\|\neq 0$ for $t>T_0$, we have $\frac{\mathrm{d}}{\mathrm{d}t}\|x_{i_0}(t)-x_{j_0}(t)\|
\leq\frac{x_{i_0}(t)-x_{j_0}(t)}{\|x_{i_0}(t)-x_{j_0}(t)\|}2q_{i_0j_0}(t)\frac{x_{j_0}(t)-x_{i_0}(t)}{\|v_{i_0}(t)-v_{j_0}(t)\|}+2nm\gamma_q\leq -2nm\gamma_q$.
When there exist at least an agent  ${i}$ such that ${i}\in N_{\tilde{i}}(t)$ and $\|x_{\tilde{i}}(t)-x_{i}(t)\|\neq0$ for $\tilde{i}\in \{i_0,j_0\}$ and either $(i_0,j_0)\notin \mathcal{G}(t)$ or $\|x_{i_0}(t)-x_{j_0}(t)\|=0$ holds, we have $\frac{\mathrm{d}}{\mathrm{d}t}\|v_{i_0}(t)-v_{j_0}(t)\|\leq 2nm\gamma_q$ for $t>T_0$. Let $T_0<t_{k3}<t_{k4}<t_{k+1,3}<t_{k+1,4}$ denote the contiguous switching times for all $k\in \{1,2,\cdots\}$ such that the case holds when $(i_0,j_0)\in \mathcal{G}(t)$ and $\|x_{i_0}(t)-x_{j_0}(t)\|\neq 0$ for all $t\in[t_{k3},t_{k4})$ and  the case when there exist at least an agent  ${i}$ such that ${i}\in N_{\tilde{i}}(t)$ and $\|x_{\tilde{i}}(t)-x_{i}(t)\|\neq0$ for $\tilde{i}\in \{i_0,j_0\}$ and either $(i_0,j_0)\notin \mathcal{G}(t)$ or $\|x_{i_0}(t)-x_{j_0}(t)\|=0$ holds for all $i,j\in \mathcal{I}$ and all $t\in[t_{k4},t_{k+1,3})$.
Note that $\|x_{i_0}(t)-x_{j_0}(t)\|<2\rho$.
Calculating $\|x_{i_0}(t)-x_{j_0}(t)\|$ based on the Newton's Law, we have that \begin{eqnarray}\label{tplink6}\begin{array}{lll}0\leq\|x_{i_0}(+\infty)-x_{j_0}(+\infty)\|\\\leq\|x_{i_0}(T_0)-x_{j_0}(T_0)\|+\sum_{k=1}^{+\infty}2nm\gamma_q(t_{k+1,3}-t_{k4})^2 \\-\sum_{k=1}^{+\infty}nm\gamma_q(t_{k4}-t_{k3})^2\\\leq 2\rho+2nm\gamma_q [\sum_{k=1}^{+\infty}(t_{k+1,3}-t_{k4})]^2-0.25nm\gamma_q[\sum_{k=1}^{+\infty}(t_{k4}-t_{k3})]^2.\end{array}\end{eqnarray} Since $\lim_{t\rightarrow+\infty}q_{i_0j_0}(t)=+\infty$, from the dynamics of $q_{ij}(t)$, we have $\sum_{k=1}^{+\infty}(t_{k4}-t_{k3})=+\infty$ and hence from (\ref{tplink6}) we have $\sum_{k=1}^{+\infty}(t_{k+1,3}-t_{k4})=+\infty$. That is, there exist a pair of agents  $i_1\neq j_1$ such that $(i_1,j_1)\notin\{(i_0,j_0),(j_0,i_0)\}$, $i_1\in \{i_0,j_0\}$ and $\lim_{t\rightarrow+\infty}q_{i_1j_1}(t)=+\infty$. Similarly, it can be proved that there exist a pair of agents $i_2\neq j_2$ such that $(i_2,j_2)\notin\{(i_0,j_0),(j_0,i_0),(i_1,j_1),(j_1,i_1)\}$, $i_2\in \{i_0,j_0,i_1\}$ and $\lim_{t\rightarrow+\infty}q_{i_2j_2}(t)=+\infty$. By analogy, it can be proved that $\lim_{t\rightarrow+\infty}q_{ij}(t)=+\infty$ for all $i,j$. Then there is a constant $T_1>0$ such that $q_{ij}(t)$ is far larger than $\rho$ for all $i,j$ and all $t>T_1$.
 Since $\bar{v}_i(t)\leq \rho$ for all $i,t$, it follows from (\ref{equationsec3}) that $\|\sum_{j\in N_i(t)}\frac{q_{ij}(t)[x_j(t)-x_i(t)]}{\|x_j(t)-x_i(t)\|}\|$ is far smaller than $\min_{j\in N_i(t)}q_{ij}(t)$ for $t>T_1$ and all $i$. Adopting a group of agents $E=\{i_1,\cdots,i_q\}$ such that $x_i(t)\in \overline{co}\{x_{i_1}(t),\cdots,x_{i_q}(t)\}$ and $x_j(t)\notin \overline{co}\{x_k(t)\mid k\in \mathcal{I}, k\neq j\}$ for all $i$ and $j\in E$ where $\overline{co}$ denotes the operator of the convex closure. It is clear that $\frac{[x_k(t)-x_{i_0}(t)]^T[x_j(t)-x_{i_0}(t)]}{\|x_j(t)-x_{i_0}(t)\|\|x_k(t)-x_{i_0}(t)\|}\geq0$ for all $i_0\in E$. If $x_j(t)\neq x_{i_0}(t)$ for some $j\in N_{i_0}(t)$, we have $\|\sum_{j\in N_i(t)}\frac{q_{ij}(t)[x_j(t)-x_i(t)]}{\|x_j(t)-x_i(t)\|}\|\geq \min_{j\in N_i(t)}q_{ij}(t)$. This yields a contradiction. Thus, $\sum_{k=1}^{+\infty}(t_{k2}-t_{k1})<+\infty$.

 Based on the above analysis, using a approach similar to the proof of Theorem \ref{theorem1}, the
team objective function (\ref{gel16}) is minimized as $t\rightarrow+\infty$.\endproof

\section{Simulations}
\begin{figure}
\centering
\includegraphics[width=3.5in]{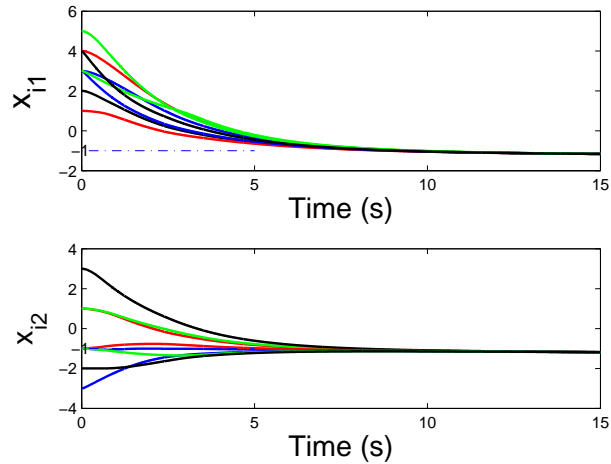} \\
\vspace{-0.2cm} \caption{State trajectories of all agents using (\ref{equationsec2}) with $p_i=1$ for all $i$.}
\label{fig:7}
\end{figure}

 Consider a multi-agent system consisting of
8 agents in a plane. The communication graph is switched among the connected subgraphs of the graph in Fig. \ref{fig:1}. The local objective
functions are  $f_1(x_1)=\frac{1}{2}x_{11}^2+\frac{1}{2}x_{12}^2,$ $f_2(x_2)=\frac{1}{2}(x_{21}+2)^2+\frac{1}{2}x_{22}^2,$
$f_3(x_3)=\frac{1}{2}x_{31}^2+\frac{1}{2}(x_{32}+2)^2,$
$f_4(x_4)=\frac{1}{2}(x_{41}+2)^2+\frac{1}{2}(x_{42}+2)^2$, $f_5(x_5)=\frac{1}{4}x_{51}^4+\frac{1}{4}x_{52}^4,$
$f_6(x_6)=\frac{1}{4}(x_{61}+2)^4+\frac{1}{4}x_{62}^4,$
$f_7(x_7)=\frac{1}{4}x_{71}^4+\frac{1}{4}(x_{72}+2)^4$ and $f_8(x_8)=\frac{1}{4}(x_{81}+2)^4+\frac{1}{4}(x_{82}+2)^4$, where $x_{i1}$ and $x_{i2}$ denote the two components of $x_i$. By simple calculations, when $\sum_{i=1}^n\nabla f_i(s)=0$, we have that $s=[-1,-1]^T$.
From Lemma \ref{lemma35},  the minimum set of the team objective function (\ref{gel16}) is $s=[-1,-1]^T$.
The simulation
results are shown in Figs. \ref{fig:6} and \ref{fig:7}. It can be observed that the team objective function (\ref{gel16}) is minimized as $t\rightarrow+\infty$, which are consistent with Theorems 1 and 2.

\begin{figure}
\begin{center}
\xymatrix{
    *++[o][F-]{1}\ar[r]\ar[d]
    & *++[o][F-]{2}\ar[r]\ar[l]
    & *++[o][F-]{3}\ar[r]\ar[l]
    & *++[o][F-]{4}\ar[d]\ar[l]
\\
    *++[o][F-]{8}\ar[u]\ar[r]
    & *++[o][F-]{7}\ar[l]\ar[r]
    & *++[o][F-]{6}\ar[l]\ar[r]
    & *++[o][F-]{5}\ar[l]\ar[u]
  }
\end{center}
\vspace{0.2cm}
\caption{One undirected graph.} \label{fig:1}
\end{figure}
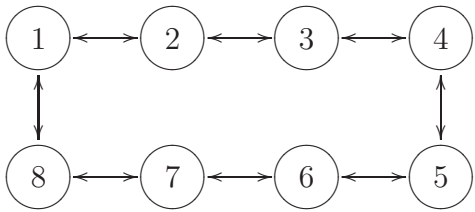
\begin{figure}
\centering
\includegraphics[width=3.5in]{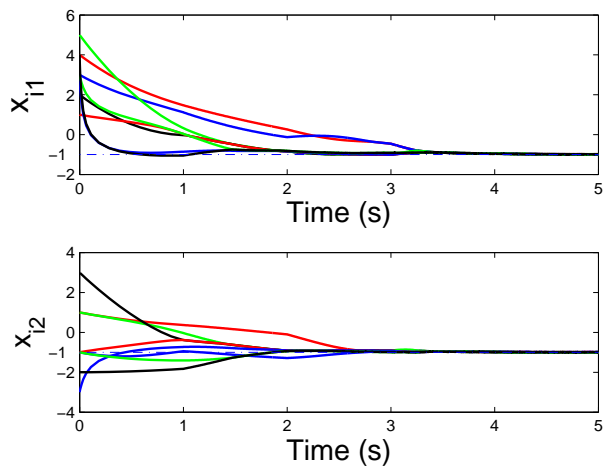} \\
\vspace{-0.2cm} \caption{State trajectories of all agents using (\ref{gel2}).}
\label{fig:6}
\end{figure}

\section{Conclusions}
In this paper, a distributed optimization problem with general differentiable convex objective functions was
studied for single-integrator and double-integrator multi-agent systems. Two distributed adaptive optimization algorithm was introduced by using the
relative information to construct the gain of the interaction term.  The analysis was performed based on the Lyapunov functions, the analysis of the system solution and the convexity of the local objective functions.  It was shown that if the gradients of the convex objective functions are continuous, the team convex objective function can be minimized as time evolves for both single-integrator and double-integrator multi-agent systems.

\end{document}